# Large deviations for stochastic flows of diffeomorphisms

AMARJIT BUDHIRAJA[1], PAUL DUPUIS[2] and VASILEIOS MAROULAS[3]

[1]*Department of Statistics and Operations Research, University of North Carolina, Chapel Hill, NC 27599, USA. E-mail: budhiraj@email.unc.edu*

[2]*Division of Applied Mathematics, Brown University, Providence, RI 02912, USA. E-mail: dupuis@dam.brown.edu*

[3]*Institute for Mathematics and its Applications, University of Minnesota, Minneapolis, MN 55455, USA. E-mail: maroulas@email.unc.edu*

A large deviation principle is established for a general class of stochastic flows in the small noise limit. This result is then applied to a Bayesian formulation of an image matching problem, and an approximate maximum likelihood property is shown for the solution of an optimization problem involving the large deviations rate function.

*Keywords:* deformable templates; diffeomorphisms; image matching; infinite-dimensional Brownian motion; infinite-dimensional SDEs; large deviations; semimartingales with a spatial parameter; small noise asymptotics; stochastic flows

## 1. Introduction

Stochastic flows of diffeomorphisms have been a subject of much research [4, 6, 13, 15]. In this paper, we are interested in an important subclass of such flows, namely the Brownian flows of diffeomorphisms (cf. [15]). Our goal is to study small noise asymptotics, specifically, the large deviation principle (LDP) for such flows.

Elementary examples of Brownian flows are those constructed by solving finite dimensional Itô stochastic differential equations. More precisely, suppose $b, f_i, i = 1, \ldots, m$, are functions from $\mathbb{R}^d \times [0,T]$ to $\mathbb{R}^d$ that are continuous in $(x,t)$ and $(k+1)$-times continuously differentiable (with uniformly bounded derivatives) in $x$. Let $\beta_1, \ldots, \beta_m$ be independent standard real Brownian motions on some filtered probability space $(\Omega, \mathcal{F}, \mathbb{P}, \{\mathcal{F}_t\})$. Then for each $s \in [0,T]$ and $x \in \mathbb{R}^d$, there is a unique continuous $\{\mathcal{F}_t\}$-adapted, $\mathbb{R}^d$-valued process $\phi_{s,t}(x)$, $s \le t \le T$, satisfying

$$\phi_{s,t}(x) = x + \int_s^t b(\phi_{s,r}(x), r) \, \mathrm{d}r + \sum_{i=1}^m \int_s^t f_i(\phi_{s,r}(x), r) \, \mathrm{d}\beta_i(r). \tag{1.1}$$

By choosing a suitable modification, $\{\phi_{s,t}, 0 \leq s \leq t \leq T\}$ defines a Brownian flow of $\mathbb{C}^k$-diffeomorphisms (see Section 2). In particular, denoting by $G^k$ the topological group of $\mathbb{C}^k$-diffeomorphisms (see Section 3 for precise definitions of the topology and the metric on $G^k$), one has that $\phi \equiv \{\phi_{0,t}, 0 \leq t \leq T\}$ is a random variable with values in the Polish space $\hat{W}_k = C([0,T]:G^k)$. For $\varepsilon \in (0,\infty)$, when $f_i$ is replaced by $\varepsilon f_i$ in (1.1), we write the corresponding flow as $\phi^\varepsilon$. Large deviations for $\phi^\varepsilon$ in $\hat{W}_k$, as $\varepsilon \to 0$, have been studied for the case $k=0$ in [3, 18] and for general $k$ in [5].

As is well known (cf. [4, 15, 16]), not all Brownian flows can be expressed as in (1.1) and in general one needs infinitely many Brownian motions to obtain a stochastic differential equation (SDE) representation for the flow. Indeed typical space–time stochastic models with a realistic correlation structure in the spatial parameter naturally lead to a formulation with infinitely many Brownian motions. One such example is given in Section 5. Thus, following Kunita's [15] notation for stochastic integration with respect to semimartingales with a spatial parameter, the study of general Brownian flows of $\mathbb{C}^k$-diffeomorphisms leads to SDEs of the form

$$d\phi_{s,t}(x) = F(\phi_{s,t}(x), dt), \qquad \phi_{s,s}(x) = x, \qquad 0 \leq s \leq t \leq T, x \in \mathbb{R}^d, \qquad (1.2)$$

where $F(x,t)$ is a $\mathbb{C}^{k+1}$-Brownian motion (see Definition 2.2). Note that such an $F$ can be regarded as a random variable with values in the Polish space $W_k = C([0,T]:\mathbb{C}^{k+1}(\mathbb{R}^d))$, where $\mathbb{C}^{k+1}(\mathbb{R}^d)$ is the space of $(k+1)$-times continuously differentiable functions from $\mathbb{R}^d$ to $\mathbb{R}^d$. Representations of such Brownian motions in terms of infinitely many independent standard real Brownian motions is well known (see, e.g., Kunita [15], Exercise 3.2.10). Indeed, one can represent $F$ as

$$F(x,t) \doteq \int_0^t b(x,r)\,dr + \sum_{i=1}^\infty \int_0^t f_i(x,r)\,d\beta_i(r), \qquad (x,t) \in \mathbb{R}^d \times [0,T], \qquad (1.3)$$

where $\{\beta_i\}_{i=1}^\infty$ is an infinite sequence of i.i.d. real Brownian motions and $b, f_i$ are suitable functions from $\mathbb{R}^d \times [0,T]$ to $\mathbb{R}^d$ (see below Definition 2.2 for details).

Letting $a(x,y,t) = \sum_{i=1}^\infty f_i(x,t)f'_i(y,t)$ for $x,y \in \mathbb{R}^d$, $t \in [0,T]$, the functions $(a,b)$ are referred to as the *local characteristics* of the Brownian motion $F$. When equation (1.2) is driven by the Brownian motion $F^\varepsilon$ with local characteristics $(\varepsilon a, b)$, we will denote the corresponding solution by $\phi^\varepsilon$. In this work we establish a large deviation principle for $(\phi^\varepsilon, F^\varepsilon)$ in $\hat{W}_{k-1} \times W_{k-1}$. Note that the LDP is established in a larger space than the one in which $(\phi^\varepsilon, F^\varepsilon)$ take values (namely, $\hat{W}_k \times W_k$). This is consistent with results in [3, 5, 18], which consider stochastic flows driven by only finitely many real Brownian motions. The main technical difficulty in establishing the LDP in $\hat{W}_k \times W_k$ is the proof of a result analogous to Proposition 4.10, which establishes tightness of certain controlled processes, when $k-1$ is replaced by $k$.

As noted above, the stochastic dynamical systems considered in this work are driven by an infinite dimensional Brownian motion. A broadly applicable approach to the study of large deviations for such systems, based on variational representations for functionals of infinite-dimensional Brownian motions, has been developed in [7, 8]. Several authors have adopted this approach to analyze the large deviation properties of a variety of models, including stochastic PDEs with random dynamic boundary conditions [22], stochastic Navier–Stokes equations [21] and infinite-dimensional SDEs with non-Lipschitz coefficients [19, 20]. The approach is a particularly attractive alternative to standard discretization/approximation methods (cf. [2]) when the state spaces are non-standard function spaces, such as the space of diffeomorphisms used in the present paper.

The proof of our main result (Theorem 3.2) proceeds by verification of a general sufficient condition obtained in [8] (see Assumption 2 and Theorem 6 therein; see also Theorem 3.6 of the current paper). The verification of this condition essentially translates into establishing weak convergence of certain stochastic flows defined via controlled analogues of the original model (see Theorem 3.5). These weak convergence proofs proceed by first establishing convergence for $N$-point motions of the flow and then using Sobolev and Rellich–Kondrachov embedding theorems (see the proof of Proposition 4.10) to argue tightness and convergence as flows. The key point here is that the estimates needed in the proofs are precisely those that have been developed in [15] for general qualitative analysis (e.g., existence, uniqueness) of the uncontrolled versions of the flows. Unlike in [3, 5, 18] (which consider only finite-dimensional flows), the proof of the LDP does not require any exponential probability estimates or discretization/approximation of the original model.

In Section 5 of this paper we study an application of these results to a problem in image analysis. Stochastic diffeomorphic flows have been suggested for modeling prior statistical distributions on the space of possible images/targets of interest in the study of nonlinear inverse problems in this field (see [12] and references therein). Along with a data model, noise corrupted observations with such a prior distribution can then be used to compute a posterior distribution on this space, the "mode" of which yields an estimate of the true image underlying the observations. Motivated by such a Bayesian procedure, a variational approach to this image matching problem has been suggested and analyzed in [12]. A goal of the current paper is to develop a rigorous asymptotic theory that relates standard stochastic Bayesian formulations of this problem, in the small noise limit, with the deterministic variational approach taken in [12]. This is done in Theorem 5.1 of Section 5.

Large deviations for finite-dimensional stochastic flows were studied for asymptotic analysis of small noise finite-dimensional anticipative SDEs in [18] and of finite-dimensional diffusions generated by $\varepsilon L_1 + L_2$, where $L_1, L_2$ are two second-order differential operators, in [5]. Analogous problems for infinite-dimensional models can be treated using the large deviation principle established in the current paper.

We now give an outline of the paper. Section 2 contains some background definitions of $\mathbb{C}^k$-Brownian motions and Brownian flows. Section 3 presents the main large deviation result of the paper. The key weak convergence needed to prove this result, Theorem 3.5, is given in Section 4. Finally, Section 5 introduces the image analysis problem and uses the results of Section 3 to obtain an asymptotic result relating the Bayesian formulation of the problem with the deterministic variational approach of [12].

We generally follow the notation of [15]. A list of standard notational conventions is given at the end of the paper, but specialized notation is as follows:

- Let $\alpha = (\alpha_1, \alpha_2, \ldots, \alpha_d)$ be a multiindex of non-negative integers and $|\alpha| = \alpha_1 + \alpha_2 + \cdots + \alpha_d$. For an $|\alpha|$-times differentiable function $f : \mathbb{R}^d \to \mathbb{R}$, set $\partial^\alpha f \doteq \partial_x^\alpha f = \frac{\partial^{|\alpha|} f}{(\partial x_1)^{\alpha_1} \cdots (\partial x_d)^{\alpha_d}}$. For such an $f$, we write $\frac{\partial f(x)}{\partial x_i}$ as $\partial_i f$. If $f \equiv (f_1, f_2, \ldots, f_d)'$ is an $|\alpha|$-times differentiable function from $\mathbb{R}^d$ to $\mathbb{R}^d$, we write $\partial^\alpha f \doteq (\partial^\alpha f_1, \partial^\alpha f_2, \ldots, \partial^\alpha f_d)'$. By convention $\partial^0 f = f$.
- For $m \geq 0$ denote by $\mathbb{C}^m$ the space of $m$-times continuously differentiable functions from $\mathbb{R}^d$ to $\mathbb{R}$.
- For any subset $A \subset \mathbb{R}^d$, $m \geq 0$ and $f \in \mathbb{C}^m$, let

$$\|f\|_{m;A} \doteq \sum_{0 \leq |\alpha| \leq m} \sup_{x \in A} |\partial^\alpha f(x)|.$$

The space $\mathbb{C}^m$ is a Fréchet space with the countable collection of seminorms $\|f\|_{m;A_n}$, $A_n = \{x : |x| \leq n\}$. In particular, it is a Polish space with a topology that corresponds to "uniform convergence on compacts".

- For $0 < \delta \leq 1$, let

$$\|f\|_{m,\delta;A} \doteq \|f\|_{m;A} + \sum_{|\alpha|=m} \sup_{x,y \in A; x \neq y} \frac{|\partial^\alpha f(x) - \partial^\alpha f(y)|}{|x-y|^\delta},$$

and

$$\mathbb{C}^{m,\delta} \doteq \{f \in \mathbb{C}^m : \|f\|_{m,\delta;A_n} < \infty \text{ for any } n \in \mathbb{N}\}.$$

The seminorms $\{\|\cdot\|_{m,\delta;A_n}, n \in \mathbb{N}\}$ make $\mathbb{C}^{m,\delta}$ a Fréchet space.

- For $m \geq 0$ denote by $\widetilde{\mathbb{C}}^m$ the space of functions $g : \mathbb{R}^d \times \mathbb{R}^d \to \mathbb{R}$ such that $g(x,y)$, $x, y \in \mathbb{R}^d$ is $m$-times continuously differentiable with respect to both $x$ and $y$. Endowed with the seminorms

$$\|g\|_{m;A_n}^\sim \doteq \sum_{0 \leq |\alpha| \leq m} \sup_{x,y \in A_n} |\partial_x^\alpha \partial_y^\alpha g(x,y)|,$$

where $n \in \mathbb{N}$, $\widetilde{\mathbb{C}}^m$ is a Fréchet space. Also, for $0 < \delta \leq 1$ let

$$\|g\|_{m,\delta;A_n}^\sim \doteq \|g\|_{m;A_n}^\sim + \sum_{|\alpha|=m} \sup_{\substack{x \neq x', y \neq y' \\ x,y,x',y' \in A_n}} \frac{|\Delta_{x,x'}^\alpha g(y) - \Delta_{x,x'}^\alpha g(y')|}{|x-x'|^\delta |y-y'|^\delta},$$

where $\Delta_{x,x'}^\alpha g(y) \doteq \hat{\partial}_{x,y}^\alpha g(x,y) - \hat{\partial}_{x',y}^\alpha g(x',y)$, $\hat{\partial}_{x,y}^\alpha g(x,y) \doteq \partial_x^\alpha \partial_y^\alpha g(x,y)$. Then

$$\widetilde{\mathbb{C}}^{m,\delta} \doteq \{g \in \widetilde{\mathbb{C}}^m ; \|g\|_{m,\delta;A_n}^\sim < \infty, \text{ for any } n \in \mathbb{N}\}$$

is a Fréchet space with respect to the seminorms $\{\|\cdot\|_{m,\delta;A_n}^\sim, n \in \mathbb{N}\}$.

- We write $\|f\|_{m;\mathbb{R}^d}$ as $\|f\|_m$. The norms $\|\cdot\|_{m,\delta}, \|\cdot\|_m^\sim, \|\cdot\|_{m,\delta}^\sim$ are to be interpreted in a similar manner.

- Let $\mathbb{C}^m(\mathbb{R}^d) \doteq \{f = (f_1, f_2, \ldots, f_d)' : f_i \in \mathbb{C}^m, i = 1, 2, \ldots, d\}$ and $\|f\|_m = \sum_{i=1}^d \|f_i\|_m$. The spaces $\mathbb{C}^{m,\delta}(\mathbb{R}^d), \widetilde{\mathbb{C}}^m(\mathbb{R}^{d \times d}), \widetilde{\mathbb{C}}^{m,\delta}(\mathbb{R}^{d \times d})$ and their corresponding norms are defined similarly. In particular, note that $h \in \widetilde{\mathbb{C}}^{m,\delta}(\mathbb{R}^{d \times d})$ is a map from $\mathbb{R}^d \times \mathbb{R}^d$ to $\mathbb{R}^{d \times d}$.

- Let $\mathbb{C}_T^{m,\delta}(\mathbb{R}^d)$ and $\widetilde{\mathbb{C}}_T^{m,\delta}(\mathbb{R}^{d \times d})$ be the classes of measurable functions $b : [0, T] \to \mathbb{C}^{m,\delta}(\mathbb{R}^d)$ and $a : [0, T] \to \widetilde{\mathbb{C}}^{m,\delta}(\mathbb{R}^{d \times d})$, respectively, such that

$$\|b\|_{T,m,\delta} \doteq \sup_{0 \leq t \leq T} \|b(t)\|_{m,\delta} < \infty \quad \text{and}$$

$$\|a\|_{T,m,\delta}^{\sim} \doteq \sup_{0 \leq t \leq T} \|a(t)\|_{m,\delta}^{\sim} < \infty.$$

- We denote the Hilbert space

$$l_2 \doteq \left\{ (x_1, x_2, \ldots) : x_i \in \mathbb{R}, i \geq 1 \text{ and } \sum_{i=1}^\infty x_i^2 < \infty \right\},$$

where the inner product on $l_2$ is defined as $\langle x, y \rangle_{l_2} \doteq \sum_{i=1}^\infty x_i y_i, x, y \in l_2$. We denote the corresponding norm as $\|\cdot\|_{l_2}$.

- Given a filtered probability space $(\Omega, \mathcal{F}, \mathbb{P}, \{\mathcal{F}_t\})$, we define

$$\mathcal{A}[l_2] \doteq \left\{ \phi \equiv \{\phi_i\}_{i=1}^\infty \ \bigg| \ \phi_i : [0, T] \to \mathbb{R} \text{ is } \{\mathcal{F}_t\}\text{-predictable for all } i \right.$$

$$\left. \text{and } \mathbb{P}\left\{ \int_0^T \|\phi(s)\|_{l_2}^2 \, ds < \infty \right\} = 1 \right\}.$$

- Define

$$S_N[l_2] \doteq \left\{ \phi \equiv \{\phi_i\}_{i=1}^\infty \in L^2([0, T] : l^2) \text{ s.t. } \int_0^T \|\phi(s)\|_{l_2}^2 \, ds \leq N \right\}.$$

Endowed with the weak topology for the Hilbert space $L^2([0, T] : l^2)$, $S_N[l_2]$ is a compact Polish space.

- Define

$$\mathcal{A}_N[l_2] \doteq \{u \in \mathcal{A}[l_2] : u(\omega) \in S_N, \mathbb{P}\text{-a.s.}\}.$$

We will always consider $S_N[l_2]$ with the weak topology when referring to convergence in distribution of $S_N[l_2]$-valued random variables.

## 2. Preliminaries

Let $\circ$ denote the composition of maps and let id denote the identity map on $\mathbb{R}^d$.

***Definition 2.1 (Stochastic flows of homeomorphisms/diffeomorphisms).*** *A collection $\{\phi_{s,t}(x): 0 \leq s \leq t \leq T, x \in \mathbb{R}^d\}$ of $\mathbb{R}^d$-valued random variables on some filtered probability space $(\Omega, \mathcal{F}, \mathbb{P}, \{\mathcal{F}_t\})$ is called a forward stochastic flow of homeomorphisms if there exists $N \in \mathcal{F}$ with $\mathbb{P}(N) = 0$ such that for any $\omega \in N^c$:*

1. *$(s,t,x) \mapsto \phi_{s,t}(x, \omega)$ is a continuous map,*
2. *$\phi_{s,u}(\omega) = \phi_{t,u}(\omega) \circ \phi_{s,t}(\omega)$ holds for all $s, t, u$, $0 \leq s \leq t \leq u \leq T$,*
3. *$\phi_{s,s}(\omega) = \mathrm{id}$ for all $s$, $0 \leq s \leq T$,*
4. *the map $\phi_{s,t}(\omega) : \mathbb{R}^d \to \mathbb{R}^d$ is an onto homeomorphism for all $s, t, 0 \leq s \leq t \leq T$.*

*If in addition $\phi_{s,t}(x, \omega)$ is $k$-times differentiable with respect to $x$ for all $s \leq t$ and the derivatives are continuous in $(s, t, x)$, it is called a stochastic flow of $\mathbb{C}^k$-diffeomorphisms.*

We now introduce a Brownian motion with a spatial parameter, with local characteristics $(a, b)$. Throughout Sections 2–4, we will assume that $(a, b) \in \widetilde{\mathbb{C}}_T^{k,\delta}(\mathbb{R}^{d \times d}) \times \mathbb{C}_T^{k,\delta}(\mathbb{R}^d)$, for some $k \in \mathbb{N}$ and $\delta \in (0, 1]$. Fix $\nu$ such that $0 < \nu < \delta$.

***Definition 2.2 ($\mathbb{C}^{k,\nu}$-Brownian motion).*** *A continuous stochastic process $\{F(t)\}_{t \geq 0}$ on some filtered probability space $(\Omega, \mathcal{F}, \mathbb{P}, \{\mathcal{F}_t\})$ with values in $\mathbb{C}^{k,\nu}(\mathbb{R}^d)$ is said to be a $\mathbb{C}^{k,\nu}$-Brownian motion with local characteristics $(a, b)$ if $F(0), F(t_{i+1}) - F(t_i), i = 0, 1, \ldots, n-1$, are independent $\mathbb{C}^{k,\nu}(\mathbb{R}^d)$-valued random variables whenever $0 \leq t_0 < t_1 < \cdots < t_n \leq T$ and if for each $x \in \mathbb{R}^d$, $M(x, t) \doteq F(x, t) - \int_0^t b(x, r)\, \mathrm{d}r$ is a continuous (d-dimensional) martingale such that $\langle\langle M(x, \cdot), M(y, \cdot) \rangle\rangle_t = \int_0^t a(x, y, r)\, \mathrm{d}r$ for all $(x, y) \in \mathbb{R}^d \times \mathbb{R}^d$.*

The existence of a $\mathbb{C}^{k,\nu}$-Brownian motion with local characteristics $(a, b)$ follows from [15] (see, e.g., Theorem 3.1.2 and Exercise 3.2.10). Indeed, for any $\gamma < \delta$ one can represent $F$ as in (1.3), where $f_i : \mathbb{R}^d \times [0, T] \to \mathbb{R}^d$ are such that for each $t \in [0, T], f_i(\cdot, t) \in \mathbb{C}^{k,\gamma}(\mathbb{R}^d)$,

$$a(x, y, t) = \sum_{i=1}^{\infty} f_i(x, t) f_i'(y, t), \qquad \text{a.e. } t,$$

and

$$\int_0^T \sum_{i=1}^{\infty} |f_i(x, r)|^2\, \mathrm{d}r \leq T \|a\|_{T,k,\delta}^{\sim} < \infty.$$

In particular, note that if $F$ is a $\mathbb{C}^{k,\nu}$-valued Brownian motion, its finite-dimensional restriction $(F(x_1, \cdot), F(x_2, \cdot), \ldots, F(x_n, \cdot))'$ is an $nd$-dimensional Brownian motion (with suitable mean and covariance) for any $(x_1, \ldots, x_n) \in \mathbb{R}^{nd}$. If $F$ is as defined by (1.3) and $\{\phi_t\}_{0 \leq t \leq 1}$ is a continuous $\mathbb{R}^d$-valued $\{\mathcal{F}_t\}$-adapted stochastic process, the stochastic integral $\int_0^t F(\phi_r, \mathrm{d}r)$ is a well-defined $d$-dimensional continuous $\{\mathcal{F}_t\}$-adapted stochastic process (see [15], Chapter 3, Section 2, pages 71–86).

***Definition 2.3.*** *Let $F$ be as in Definition 2.2. Then for each $s \in [0, T]$ and $x \in \mathbb{R}^d$, there is a unique continuous $\mathcal{F}_t$-adapted, $\mathbb{R}^d$-valued process $\phi_{s,t}(x), s \leq t \leq T$ satisfying* $\phi_{s,t}(x) = x + \int_s^t F(\phi_{s,r}(x), \mathrm{d}r), t \in [s,T]$. *This stochastic process is called the solution of Itô's stochastic differential equation based on the Brownian motion $F$.*

From [15], Theorem 4.6.5, it follows that $\{\phi_{s,t}\}_{0 \le s \le t \le T}$ as introduced in Definition 2.3 has a modification that is a forward stochastic flow of $\mathbb{C}^k$-diffeomorphisms.

## 3. Large deviation principle

Given $\varepsilon > 0$, let $F^\varepsilon$ be a $\mathbb{C}^{k,\nu}$-Brownian motion on some filtered probability space $(\Omega, \mathcal{F}, \mathbb{P}, \{\mathcal{F}_t\})$, with local characteristics $(\varepsilon a, b)$, where $(k, \nu)$ and $(a, b)$ are as in Section 2. Without loss of generality we assume that $F^\varepsilon$ is represented as

$$F^\varepsilon(x,t) \doteq \int_0^t b(x,r)\,\mathrm{d}r + \sqrt{\varepsilon} \sum_{l=1}^\infty \int_0^t f_l(x,r)\,\mathrm{d}\beta_l(r), \qquad (x,t) \in \mathbb{R}^d \times [0,T], \qquad (3.1)$$

where $(\beta_l, f_l)_{l \ge 1}$ are as in Section 2. Note in particular that

$$F^\varepsilon(x,t) - \int_0^t b(x,r)\,\mathrm{d}r = \sqrt{\varepsilon} M(x,t).$$

With an abuse of notation, when $\varepsilon = \varepsilon_n$ we write $F^\varepsilon$ as $F^n$. Observe that $\langle\!\langle M(x,\cdot), \beta_l(\cdot) \rangle\!\rangle_t = \int_0^t f_l(x,r)\,\mathrm{d}r$ for all $t \in [0,T]$, a.s. Let $\phi^\varepsilon \equiv \{\phi^\varepsilon_{s,t}(x), 0 \le s \le t \le T, x \in \mathbb{R}^d\}$ be the forward stochastic flow of $\mathbb{C}^k$-diffeomorphisms based on $F^\varepsilon$. With another abuse of notation, we write $\phi^\varepsilon_{0,t}$ as $\phi^\varepsilon_t$ and $\phi^\varepsilon = \{\phi^\varepsilon_t(x), 0 \le t \le T, x \in \mathbb{R}^d\}$.

The goal of this paper is to show that the family $(\phi^\varepsilon, F^\varepsilon)_{\varepsilon > 0}$ satisfies an LDP on a suitable function space, as $\varepsilon \to 0$. We begin by recalling the definition of a rate function.

**Definition 3.1.** *Let $\mathcal{E}$ be a Polish space. A function $I : \mathcal{E} \to [0, \infty]$ is called a rate function on $\mathcal{E}$, if for each $M < \infty$ the level set $\{x \in \mathcal{E} : I(x) \le M\}$ is a compact subset of $\mathcal{E}$.*

We remark that some authors define a rate function by replacing the requirement of compactness of level sets by the requirement that these level sets be closed and refer to a rate function satisfying the property of compact level sets as a "good rate function".

For $m \in \mathbb{N}$, let $G^m$ be the group of $\mathbb{C}^m$-diffeomorphisms on $\mathbb{R}^d$. $G^m$ is endowed with the metric

$$d_m(\phi, \psi) = \lambda_m(\phi, \psi) + \lambda_m(\phi^{-1}, \psi^{-1}), \qquad (3.2)$$

where

$$\lambda_m(\phi, \psi) = \sum_{|\alpha| \le m} \rho(\partial^\alpha \phi, \partial^\alpha \psi),$$

(3.3)

$$\rho(\phi,\psi) = \sum_{N=1}^{\infty} \frac{1}{2^N} \frac{\sup_{|x|\leq N} |\phi(x) - \psi(x)|}{1 + \sup_{|x|\leq N} |\phi(x) - \psi(x)|}.$$

Under this metric $G^m$ is a Polish space. Let $\hat{W}_m \doteq C([0,T]:G^m)$ be the set of all continuous maps from $[0,T]$ to $G^m$ and $W_m \doteq C([0,T]:\mathbb{C}^m(\mathbb{R}^d))$ be the set of all continuous maps from $[0,T]$ to $\mathbb{C}^m(\mathbb{R}^d)$. The space $\hat{W}_m$ endowed with the metric $\hat{d}_m(\phi,\psi) = \sup_{0\leq t\leq T} d_m(\phi(t),\psi(t))$ and the space $W_m$ with the metric $\bar{d}_m(\phi,\psi) = \sup_{0\leq t\leq T} \lambda_m(\phi(t),\psi(t))$ are Polish spaces. Note that $(\phi^\varepsilon, F^\varepsilon)$ belongs to $\hat{W}_k \times W_k \subseteq \hat{W}_{k-1} \times W_{k-1} \subseteq W_{k-1} \times W_{k-1}$. We will show that the pair $(\phi^\varepsilon, F^\varepsilon)_{\varepsilon>0}$ satisfies LDPs in both of the spaces $\hat{W}_{k-1} \times W_{k-1}$ and $W_{k-1} \times W_{k-1}$, with a rate function $I$ that is introduced below.

Let $u \equiv \{u_l\}_{l=1}^{\infty} \in \bigcup_{N\geq 1} \mathcal{A}_N[l_2]$. Given any such control, we want to construct a corresponding controlled flow in the form of a perturbed analogue of (3.1). Observe that $Z_t \doteq \sum_{l=1}^{\infty} \int_0^t u_l(s)\,\mathrm{d}\beta_l(s)$ is a continuous square integrable martingale. For any $\gamma < \delta$ one can find $b_u:\mathbb{R}^d \times [0,T] \times \Omega \to \mathbb{R}^d$ such that $b_u(t,\omega) \in \mathbb{C}^{k,\gamma}(\mathbb{R}^d)$ for a.e. $(t,\omega)$, such that for each $x \in \mathbb{R}^d$, $b_u(x,\cdot)$ is predictable, and such that $\int_0^t b_u(x,s)\,\mathrm{d}s = \langle\!\langle Z, M(x,\cdot)\rangle\!\rangle_t$ for each $(x,t) \in \mathbb{R}^d \times [0,T]$. In particular, for each $x \in \mathbb{R}^d$, $b_u(x,t) \doteq \sum_{l=1}^{\infty} u_l(t) f_l(x,t)$ a.e. $(t,\omega)$. Furthermore, for some $c \in (0,\infty)$,

$$\|b_u(t)\|_{k,\gamma}^2 \leq c\|a\|_{T,k,\delta}^{\sim} \sum_{l=1}^{\infty} |u_l(t)|^2, \qquad [\mathrm{d}t \otimes \mathbb{P}]\text{-a.e. in }(t,\omega). \tag{3.4}$$

The proofs of these statements follow along the lines of Exercise 3.2.10 and Lemma 3.2.3 of [15]. Next, define

$$F^{0,u}(x,t) \doteq \int_0^t b_u(x,s)\,\mathrm{d}s + \int_0^t b(x,s)\,\mathrm{d}s. \tag{3.5}$$

It follows that $F^{0,u}(\cdot,t)$ is a $\mathbb{C}^{k,\gamma}(\mathbb{R}^d)$-valued continuous adapted stochastic process. Let $\hat{b}_u \doteq b_u + b$ and for $(t_0,x) \in [0,T] \times \mathbb{R}^d$ let $\{\phi_{t_0,t}^{0,u}(x)\}_{t_0\leq t\leq T}$ be the unique solution of the equation

$$\phi_{t_0,t}^{0,u}(x) \doteq x + \int_{t_0}^t \hat{b}_u(\phi_{t_0,r}^{0,u}(x),r)\,\mathrm{d}r, \qquad t \in [t_0,T]. \tag{3.6}$$

From [15], Theorem 4.6.5, it follows that $\{\phi_{s,t}^{0,u}, 0 \leq s \leq t \leq T\}$ is a forward flow of $\mathbb{C}^k$-diffeomorphisms.

For $(\phi^0, F^0) \in \hat{W}_k \times W_k$ define

$$I(\phi^0, F^0) \doteq \inf_{u \in \mathcal{L}(\phi_0, F_0)} \frac{1}{2} \int_0^T \|u(s)\|_{l_2}^2\,\mathrm{d}s, \tag{3.7}$$

where $\mathcal{L}(\phi^0, F^0) = \{u \in L^2([0,T]:l_2) \mid (\phi^0, F^0) = (\phi^{0,u}, F^{0,u})\}$. Note in particular that $u$ in (3.7) is deterministic. If $(\phi^0, F^0) \in (W_{k-1} \times W_{k-1}) \setminus (\hat{W}_k \times W_k)$ then we set $I(\phi^0, F^0) =$ $\infty$. We denote the restriction of $I$ to $\hat{W}_{k-1} \times W_{k-1}$ by the same symbol. The following is the main result of the section.

**Theorem 3.2 (Large deviation principle).** *The family $(\phi^\varepsilon, F^\varepsilon)_{\varepsilon>0}$ satisfies an LDP in the spaces $\hat{W}_{k-1} \times W_{k-1}$ and $W_{k-1} \times W_{k-1}$ with rate function $I$.*

Let $\{u_n\}_{n=1}^\infty (u_n \equiv \{u_l^n\}_{l=1}^\infty)$ be a sequence in $\mathcal{A}_N[l_2]$ for some fixed $N < \infty$. Let $\{\varepsilon_n\}_{n\geq 0}$ be a sequence such that $\varepsilon_n \geq 0$ for each $n$ and $\varepsilon_n \to 0$ as $n \to \infty$. Note that we allow $\varepsilon_n = 0$ for all $n$. Recall that $M(x,t) = \sum_{i=1}^\infty \int_0^t f_i(x,r) \, \mathrm{d}\beta_i(r), (x,t) \in \mathbb{R}^d \times [0,T]$. Define

$$\hat{F}^n(x,t) \doteq \int_0^t \hat{b}_{u_n}(x,r) \, \mathrm{d}r + \sqrt{\varepsilon_n} M(x,t) \tag{3.8}$$

and let $\phi^n$ be the solution to

$$\phi^n_t(x) = x + \int_0^t \hat{b}_{u_n}(\phi^n_r(x), r) \, \mathrm{d}r + \sqrt{\varepsilon_n} \int_0^t M(\phi^n_r(x), \mathrm{d}r). \tag{3.9}$$

Clearly $\hat{F}^n \in W_k$, and from [15], Theorem 4.6.5, equation (3.9) has a unique solution $\phi^n \in \hat{W}_k$ a.s. We next introduce some basic weak convergence definitions.

**Definition 3.3.** *Let $u \in \mathcal{A}_N[l_2]$ and $\{\phi^n\}$ be as above. Let $\hat{\mathbb{P}}^n_{k-1}, \hat{\mathbb{P}}^0_{k-1}$ be the measures induced by $(\phi^n, \hat{F}^n), (\phi^{0,u}, F^{0,u})$, respectively, on $\hat{W}_{k-1} \times W_{k-1}$. Thus for $A \in \mathcal{B}(\hat{W}_{k-1} \times W_{k-1})$,*

$$\hat{\mathbb{P}}^n_{k-1}(A) = \mathbb{P}((\phi^n, \hat{F}^n) \in A), \qquad \hat{\mathbb{P}}^0_{k-1}(A) = \mathbb{P}((\phi^{0,u}, F^{0,u}) \in A).$$

*The sequence $\{(\phi^n, \hat{F}^n)\}_{n\geq 1}$ is said to converge weakly as $G^{k-1}$-flows to $(\phi^{0,u}, F^{0,u})$ as $n \to \infty$ if $\hat{\mathbb{P}}^n_{k-1}$ converges weakly to $\hat{\mathbb{P}}^0_{k-1}$ as $n \to \infty$.*

**Definition 3.4.** *Let $\mathbb{P}^n_{k-1}, \mathbb{P}^0_{k-1}$ be the measures induced by $(\phi^n, \hat{F}^n), (\phi^{0,u}, F^{0,u})$, respectively, on $W_{k-1} \times W_{k-1}$. The sequence $\{(\phi^n, \hat{F}^n)\}_{n\geq 1}$ is said to converge weakly as $\mathbb{C}^{k-1}$-flows to $(\phi^{0,u}, F^{0,u})$ as $n \to \infty$ if $\mathbb{P}^n_{k-1}$ converges weakly to $\mathbb{P}^0_{k-1}$ as $n \to \infty$.*

As noted in the Introduction, proofs of large deviations properties based on the general framework developed in [8] essentially reduce to weak convergence questions for controlled analogues of the original process. For our problem the following theorem gives the needed result. The proof is given in the next section.

**Theorem 3.5.** *Let $\{u_n\}$ converge to $u$ in distribution as an $S_N[l_2]$-valued sequence of random variables. Then the sequence $\{(\phi^n, \hat{F}^n)\}_{n\geq 1}$ converges weakly as $\mathbb{C}^{k-1}$-flows and $G^{k-1}$-flows to the pair $(\phi^{0,u}, F^{0,u})$ as $n \to \infty$.*

The following theorem is taken from [8]. Let $\mathbb{R}^\infty$ denote the product space of countably many copies of the real line. Then $S \equiv C([0,T]:\mathbb{R}^\infty)$ (with the usual topology) is a Polish space and $\beta \equiv \{\beta_i\}_{i=1}^\infty$ is an $S$-valued random variable.

**Theorem 3.6.** *Let $\mathcal{E}$ be a Polish space, let $\{\mathcal{G}^\varepsilon\}_{\varepsilon \geq 0}$ be a collection of measurable maps from $(S, \mathcal{B}(S))$ to $(\mathcal{E}, \mathcal{B}(\mathcal{E}))$, and let $X^\varepsilon = \mathcal{G}^\varepsilon(\sqrt{\varepsilon}\beta)$. Suppose that there exists a measurable map $\mathcal{G}^0 : S \to \mathcal{E}$ such that for every $N < \infty$ the set $\Gamma_N \doteq \{\mathcal{G}^0(\int_0^\cdot u(s) \, ds) : u \in S_N[l_2]\}$ is a compact subset of $\mathcal{E}$. For $f \in \mathcal{E}$ let $\mathcal{C}_f = \{u \in L^2([0,T] : l^2) : f = \mathcal{G}^0(\int_0^\cdot u_s \, ds)\}$. Then $\hat{I}$ defined by*

$$\hat{I}(f) = \inf_{u \in \mathcal{C}_f} \left\{ \frac{1}{2} \int_0^T \|u(s)\|_{l_2}^2 \, ds \right\}, \qquad f \in \mathcal{E},$$

*is a rate function on $\mathcal{E}$. Furthermore, suppose that for all $N < \infty$ and families $\{u^\varepsilon\} \subset \mathcal{A}_N[l_2]$ such that $u^\varepsilon$ converges in distribution to some $u \in \mathcal{A}_N[l_2]$, we have that $\mathcal{G}^\varepsilon(\sqrt{\varepsilon}\beta(\cdot) + \int_0^\cdot u^\varepsilon(s) \, ds) \to \mathcal{G}^0(\int_0^\cdot u(s) \, ds)$ in distribution as $\varepsilon \to 0$. Then the family $\{X^\varepsilon, \varepsilon > 0\}$ satisfies the LDP on $\mathcal{E}$, as $\varepsilon \to 0$, with rate function $\hat{I}$.*

**Proof of Theorem 3.2.** We will only show that the sequence $(\phi^\varepsilon, F^\varepsilon)$ satisfies an LDP in $\hat{W}_{k-1} \times W_{k-1}$ with rate function $I$ defined as in (3.7). The LDP in $W_{k-1} \times W_{k-1}$ follows similarly. Let $\mathcal{G}^\varepsilon : S \to \hat{W}_{k-1} \times W_{k-1}$ be a measurable map such that $\mathcal{G}^\varepsilon(\sqrt{\varepsilon}\beta) = (\phi^\varepsilon, F^\varepsilon)$ a.s., where $F^\varepsilon$ is given by (3.1) and $\phi^\varepsilon$ is the associated flow based on $F^\varepsilon$. Define $\mathcal{G}^0 : S \to \hat{W}_{k-1} \times W_{k-1}$ by $\mathcal{G}^0(\int_0^\cdot u(s) \, ds) = (\phi^0, F^0)$ if $u \in L^2([0,T] : l_2)$ and with $\phi^0, F^0$ as defined in (3.6) and (3.5), respectively. We set $\mathcal{G}^0(f) = 0$ for all other $f \in S$.

Fix $N < \infty$ and consider $\Gamma_N = \{\mathcal{G}^0(\int_0^\cdot u(s) \, ds), u \in S_N[l_2]\}$. We first show that $\Gamma_N$ is a compact subset of $\hat{W}_{k-1} \times W_{k-1}$. It suffices to show that if $u_n, u \in S_N[l_2]$ are such that $u_n \to u$, then $\mathcal{G}^0(\int_0^\cdot u_n(s) \, ds) \to \mathcal{G}^0(\int_0^\cdot u(s) \, ds)$ in $\hat{W}_{k-1} \times W_{k-1}$. This is immediate from Theorem 3.5 on noting that $\mathcal{G}^0(\int_0^\cdot u_n(s) \, ds) = (\phi^n, \hat{F}^n)$, where $\phi^n, \hat{F}^n$ are as in (3.9) and (3.8), respectively, with $\varepsilon_n = 0$ and $\mathcal{G}^0(\int_0^\cdot u(s) \, ds) = (\phi^{0,u}, F^{0,u})$, where $\phi^{0,u}, F^{0,u}$ are as in (3.6) and (3.5), respectively.

Next let $\{u_n\} \subset \mathcal{A}_N[l_2]$ and $\varepsilon_n \in (0, \infty)$ be such that $\varepsilon_n \to 0$ and $u_n$ converges in distribution to some $u$ as $n \to \infty$. In order to complete the proof, it is enough, in view of Theorem 3.6 and the definition (3.7), to show that $\mathcal{G}^{\varepsilon_n}(\sqrt{\varepsilon_n}\beta + \int_0^\cdot u_n(s) \, ds) \to \mathcal{G}^0(\int_0^\cdot u(s) \, ds)$ in $\hat{W}_{k-1} \times W_{k-1}$, as $n \to \infty$. An application of Girsanov's theorem shows that $\mathcal{G}^{\varepsilon_n}(\sqrt{\varepsilon_n}\beta + \int_0^\cdot u_n(s) \, ds) = (\phi^n, \hat{F}^n)$, where $\phi^n, \hat{F}^n$ are defined as in (3.9) and (3.8), respectively. Also $\mathcal{G}^0(\int_0^\cdot u(s) \, ds) = (\phi^{0,u}, F^{0,u})$, where $\phi^{0,u}, F^{0,u}$ are the same as in (3.6) and (3.5), respectively. The result now follows from Theorem 3.5. □

## 4. Proof of Theorem 3.5

This section will present the proof of Theorem 3.5. It is worth recalling the assumptions that will be in effect for this section, which are that $\{u_n\}$ is converging to $u$ in distribution as an $S_N[l_2]$-valued sequence of random variables and that $(a, b) \in \widetilde{\mathbb{C}}_T^{k,\delta}(\mathbb{R}^{d \times d}) \times \mathbb{C}_T^{k,\delta}(\mathbb{R}^d)$, for some $k \in \mathbb{N}$ and $\delta \in (0, 1]$.

We begin by introducing the $(m+p)$-point motion of the flow and the related notion of "convergence as diffusions". Let $\mathbf{x} = (x_1, x_2, \ldots, x_m)$ and $\mathbf{y} = (y_1, y_2, \ldots, y_p)$ be arbitrary fixed points in $\mathbb{R}^{d\times m}$ and $\mathbb{R}^{d\times p}$, respectively. Set

$$\phi_t^n(\mathbf{x}) = (\phi_t^n(x_1), \phi_t^n(x_2), \ldots, \phi_t^n(x_m))$$

and

$$\hat{F}^n(\mathbf{y}, t) = (\hat{F}^n(y_1, t), \hat{F}^n(y_2, t), \ldots, \hat{F}^n(y_p, t)).$$

Then the pair $\{\phi_t^n(\mathbf{x}), \hat{F}^n(\mathbf{y}, t)\}$ is a continuous stochastic process with values in $\mathbb{R}^{d\times m} \times \mathbb{R}^{d\times p}$ and is called an $(m+p)$-point motion of the flow. Let $V_m \doteq C([0,T] : \mathbb{R}^{d\times m})$ be the Fréchet space of all continuous maps from $[0,T]$ to $\mathbb{R}^{d\times m}$, with the usual seminorms, and let $V_{m,p} = V_m \times V_p$ be the product space.

**Definition 4.1.** *Let $\mathbb{P}^n_{(\mathbf{x},\mathbf{y})}$ and $\mathbb{P}^0_{(\mathbf{x},\mathbf{y})}$ be the measures induced by $(\phi^n(\mathbf{x}), \hat{F}^n(\mathbf{y}))$ and $(\phi^{0,u}(\mathbf{x}), F^{0,u}(\mathbf{y}))$, respectively, on $V_{m,p}$. Thus for $A \in \mathcal{B}(V_{m,p})$,*

$$\mathbb{P}^n_{(\mathbf{x},\mathbf{y})} = \mathbb{P}((\phi^n(\mathbf{x}), \hat{F}^n(\mathbf{y})) \in A), \qquad \mathbb{P}^0_{(\mathbf{x},\mathbf{y})} = \mathbb{P}((\phi^{0,u}(\mathbf{x}), F^{0,u}(\mathbf{y})) \in A).$$

*The sequence $\{(\phi^n, \hat{F}^n)\}_{n\geq 1}$ is said to converge weakly as diffusions to $(\phi^{0,u}, F^{0,u})$ as $n \to \infty$ if $\mathbb{P}^n_{(\mathbf{x},\mathbf{y})}$ converges weakly to $\mathbb{P}^0_{(\mathbf{x},\mathbf{y})}$ as $n \to \infty$ for each $(\mathbf{x},\mathbf{y}) \in \mathbb{R}^{d\times m} \times \mathbb{R}^{d\times p}$, and $m, p = 1, 2, \ldots$.*

The following well-known result (cf. [15], Theorem 5.1.1) is a key ingredient to the proof of Theorem 3.5.

**Theorem 4.2.** *The family of probability measures $\hat{\mathbb{P}}^n_{k-1}$ (respectively, $\mathbb{P}^n_{k-1}$) converges weakly to probability measures $\hat{\mathbb{P}}^0_{k-1}$ (respectively, $\mathbb{P}^0_{k-1}$), as $n \to \infty$ if and only if the following two conditions are satisfied:*

1. *the sequence $\{(\phi^n, \hat{F}^n)\}_{n\geq 1}$ converges weakly as diffusions to $(\phi^{0,u}, F^{0,u})$ as $n \to \infty$,*
2. *the sequence $\{\hat{\mathbb{P}}^n_{k-1}\}$ (respectively, $\{\mathbb{P}^n_{k-1}\}$) is tight.*

We will show first that under the condition of Theorem 3.5 the sequence $\{(\phi^n, \hat{F}^n)\}_{n\geq 1}$ converges weakly as diffusions to $(\phi^{0,u}, F^{0,u})$ as $n \to \infty$. We begin with the following lemma.

**Lemma 4.3.** *For each $x \in \mathbb{R}^d$*

$$\mathbb{E} \sup_{0\leq t\leq T} \left| \sum_{k=1}^\infty \int_0^t f_k(x,s)\,\mathrm{d}\beta_k(s) \right|^2 < \infty, \tag{4.1}$$

$$\sup_n \mathbb{E} \sup_{0\leq t\leq T} \left| \sum_{k=1}^\infty \int_0^t f_k(\phi_s^n(x), s)\,\mathrm{d}\beta_k(s) \right|^2 < \infty. \tag{4.2}$$

**Proof.** We will only prove (4.2). The proof of (4.1) follows in a similar manner. From the Bürkholder–Davis–Gundy inequality the left-hand side of (4.2) is bounded by

$$c_1 \mathbb{E} \left| \sum_{l=1}^{\infty} \int_0^T \text{Tr}(f_l f_l')(\phi_r(x), r) \, dr \right| = c_1 \mathbb{E} \left| \int_0^T \text{Tr}(a(\phi_r(x), \phi_r(x), r)) \, dr \right|$$

$$\leq c_2 \|a\|_{T,k,\delta}^{\sim}.$$

The last expression is finite since $a$ belongs to $\tilde{\mathbb{C}}_T^{k,\delta}(\mathbb{R}^{d \times d})$. □

An immediate consequence of Lemma 4.3 is the following corollary (cf. (3.8), (3.9)).

**Corollary 4.4.** *For each $x \in \mathbb{R}^d$ and $t \in [0, T]$,*

$$\hat{F}^n(x, t) = \int_0^t \hat{b}_{u_n}(x, r) \, dr + S_n(x, t)$$

*and*

$$\phi_t^n(x) = x + \int_0^t \hat{b}_{u_n}(\phi_r^n(x), r) \, dr + T_n(x, t),$$

*where $S_n(x, \cdot)$ and $T_n(x, \cdot)$ are continuous stochastic processes with values in $\mathbb{R}^d$, satisfying $\sup_{0 \leq t \leq T}\{|S_n(x, t)| + |T_n(x, t)|\} \to 0$ in probability as $n \to \infty$.*

The following lemma, showing the tightness of $\mathbb{P}^n_{(\mathbf{x}, \mathbf{y})}$, plays an important role in the proof of the weak convergence as diffusions.

**Lemma 4.5.** *For each $x \in \mathbb{R}^d$ the sequence $\{(\phi^n(x), \hat{F}^n(x))\}_{n \geq 1}$ is tight in $C([0, T] : \mathbb{R}^d \times \mathbb{R}^d)$.*

**Proof.** We will only argue the tightness of $\{\phi^n(x)\}$. Tightness of $\{\hat{F}^n(x)\}$ is proved similarly. Corollary 4.4 yields that $T_n(x, \cdot)$ is tight in $C([0, T] : \mathbb{R}^d)$. Thus it suffices to show the tightness of $\{\int_0^{\cdot} \hat{b}_{u_n}(\phi_r^n(x), r) \, dr\}$. Fix $p > 0$. From the Cauchy–Schwarz inequality, (3.4) and recalling that $u_n \in \mathcal{A}_N[l_2]$, $\mathbb{E}| \int_s^t \hat{b}_{u_n}(\phi_r^n(x), r) \, dr|^p$ can be bounded by

$$\mathbb{E}\left[ \int_s^t |\hat{b}_{u_n}(\phi_r^n(x), r)|^2 \, dr \right]^{p/2} (t-s)^{p/2} \leq c_1 \{\|a\|_{T,k,\delta}^{\sim} + \|b\|_{T,k,\delta}^2\}^{p/2} (t-s)^{p/2}$$

$$\leq c_2 (t-s)^{p/2},$$

the result follows. □

**Proposition 4.6.** *Let $u_n \to u$ in distribution as $S_N[l_2]$-valued random variables. Then the sequence $\{(\phi^n, \hat{F}^n)\}_{n \geq 1}$ converges weakly as diffusions to $(\phi^{0,u}, F^{0,u})$ as $n \to \infty$.*

**Proof.** We claim that it suffices to show that for each $t \in [0, T]$, the map

$$(\xi, v) \mapsto \int_0^t \hat{b}_v(\xi_s, s) \, \mathrm{d}s \tag{4.3}$$

from $C([0,T] : \mathbb{R}^d) \times S_N[l_2]$ to $\mathbb{R}^d$ is continuous. To see the claim, note that in view of the tightness established in Lemma 4.5, the proposition will follow if any weak limit point $(\bar{\phi}, \bar{F}, \bar{u})$ of $(\phi^n, \hat{F}^n, u_n)$ satisfies, for each fixed $(t, x) \in [0, T] \times \mathbb{R}^d$,

$$\bar{F}(x, t) = \int_0^t \hat{b}_{\bar{u}}(x, r) \, \mathrm{d}r, \qquad \bar{\phi}_t(x) = x + \int_0^t \hat{b}_{\bar{u}}(\bar{\phi}_r(x), r) \, \mathrm{d}r, \qquad \text{a.s.} \tag{4.4}$$

Now fix a weakly convergent subsequence and $(t, x) \in [0, T] \times \mathbb{R}^d$. From (4.3) $(\phi_t^n(x), \int_0^t \hat{b}_{u_n}(\phi_r^n(x), r) \, \mathrm{d}r)$ converges weakly in $\mathbb{R}^d \times \mathbb{R}^d$ to $(\bar{\phi}_t(x), \int_0^t \hat{b}_{\bar{u}}(\bar{\phi}_r(x), r) \, \mathrm{d}r)$. The second equality in (4.4) is now an immediate consequence of the second equality in Corollary 4.4. The first equality in (4.4) is proved similarly on noting that (4.3) in particular implies that the map $v \mapsto \int_0^t \hat{b}_v(x, s) \, \mathrm{d}s$, from $S_N[l_2]$ to $\mathbb{R}^d$, is continuous.

We now prove (4.3). Let $(\xi_n, v_n) \to (\xi, v)$ in $C([0, T] : \mathbb{R}^d) \times S_N[l_2]$. Then,

$$\left| \int_0^t (\hat{b}_{v_n}(\xi_s^n, s) - \hat{b}_v(\xi_s, s)) \, \mathrm{d}s \right| \leq \left| \int_0^t (\hat{b}_{v_n}(\xi_s^n, s) - \hat{b}_{v_n}(\xi_s, s)) \, \mathrm{d}s \right|$$
$$+ \left| \int_0^t (\hat{b}_{v_n}(\xi_s, s) - \hat{b}_v(\xi_s, s)) \, \mathrm{d}s \right| \tag{4.5}$$
$$\equiv L_1 + L_2.$$

For each $x \in \mathbb{R}^d$ we have that

$$\left| \int_0^t (\hat{b}_{v_n}(x, s) - \hat{b}_v(x, s)) \, \mathrm{d}s \right| = \left| \sum_{l=1}^\infty \int_0^t f_l(x, s)(v_l^n(s) - v_l(s)) \, \mathrm{d}s \right| \to 0, \tag{4.6}$$

since $v_n \to v$ weakly in $L^2([0, T] : l_2)$ and

$$\sum_{l=1}^\infty \int_0^t |f_l(x, s)|^2 \, \mathrm{d}s \leq T \|a\|_{T, k, \delta}^{\sim} < \infty.$$

Furthermore from (3.4) (recall $k \geq 1$) we have that for some $c_1 \in (0, \infty)$ and all $x, y \in \mathbb{R}^d, 0 \leq t \leq T$,

$$\left| \int_0^t (\hat{b}_{v_n}(x, s) - \hat{b}_{v_n}(y, s)) \, \mathrm{d}s \right| \leq |x - y| \int_0^t (\|b_{v_n}(s)\|_{k, \gamma} + \|b(s)\|_{k, \gamma}) \, \mathrm{d}s$$
$$\leq c_1 |x - y|. \tag{4.7}$$

Using the Ascoli–Arzelà theorem (in the spatial variable) and equation (4.6), (4.7) now yields that the expression on the left-hand side of (4.6) converges to 0 uniformly for $x$ in compact subsets of $\mathbb{R}^d$. Thus $L_2 \to 0$ as $n \to \infty$. Following similar arguments $L_1$ is bounded by $c_2 \sup_{0 \leq s \leq T} |\xi_s^n - \xi_s|$, which converges to 0 as $n \to \infty$. Hence (4.5) converges to 0 as $n \to \infty$ and the result follows. □

We next show the tightness of the family of probability measures $\{\mathbb{P}_{k-1}^n\}$. Key ingredients in the proof are the following uniform $L^p$-estimates on $\partial^\alpha \hat{F}^n(x,t)$ and $\partial^\alpha \phi_t^n(x)$.

**Lemma 4.7.** *For each $p \geq 1$ there exists $k_1 \in (0, \infty)$ such that for all $t, t' \in [0, T], x \in \mathbb{R}^d$, $n \geq 1$ and $|\alpha| \leq k$:*

$$\mathbb{E}|\partial^\alpha \hat{F}^n(x,t) - \partial^\alpha \hat{F}^n(x,t')|^p \leq k_1 |t - t'|^{p/2}. \tag{4.8}$$

**Proof.** Fix a multiindex $\alpha$ such that $|\alpha| \leq k$ and $p \geq 1$. Using the Bürkholder–Davis–Gundy inequality for the martingale $\partial^\alpha M(x, \cdot)$ and the fact that $a \in \widetilde{\mathbb{C}}_T^{k,\delta}(\mathbb{R}^d)$, we obtain that for some $c_1 \in (0, \infty)$ and all $x \in \mathbb{R}^d, t, t' \in [0, T]$,

$$\mathbb{E}|\partial^\alpha M(x,t) - \partial^\alpha M(x,t')|^p \leq c_1 |t - t'|^{p/2}. \tag{4.9}$$

Recalling that $\hat{b}_{u_n}(\cdot, t) \in \mathbb{C}^{k,\gamma}(\mathbb{R}^d)$ a.e. $(t, \omega)$ and using (3.4) we get

$$\int_0^t \sup_{x \in \mathbb{R}^d} |\partial^\alpha \hat{b}_{u_n}(x, r)| \, dr < \infty \qquad \text{a.e.,}$$

and thus $\partial^\alpha \int_0^t b_{u_n}(x, r) \, dr = \int_0^t \partial^\alpha b_{u_n}(x, r) \, dr$ a.e. An application of the Cauchy–Schwarz inequality and (3.4) now gives, for some $c_2 \in (0, \infty)$,

$$\mathbb{E}\left|\partial^\alpha \int_{t'}^t b_{u_n}(x, r) \, dr\right|^p \leq c_2 |t - t'|^{p/2}. \tag{4.10}$$

Equation (4.8) is an immediate consequence of (4.9) and (4.10). □

For $g: \mathbb{R}^d \times [0, T] \to \mathbb{R}^d$, let $\nabla_y g(y, r)$ be the $d \times d$ matrix with entries $[\nabla_y g(y, r)]_{ij} = \frac{\partial}{\partial y_j} g_i(y, r)$. Differentiating with respect to $x_1$ in (3.9) we obtain

$$\partial_1 \phi_t^n(x) = \partial_1 x + \int_0^t [\nabla_y \hat{b}_{u_n}(\phi_r^n(x), r) \cdot \partial_1 \phi_r^n(x)] \, dr$$

$$+ \sqrt{\varepsilon_n} \int_0^t \nabla_y M(\phi_r^n(x), dr) \cdot \partial_1 \phi_r^n(x)$$

$$= \partial_1 x + \int_0^t \nabla_y \hat{F}^n(\phi_r^n(x), dr) \cdot \partial_1 \phi_r^n(x).$$

By repeated differentiation one obtains the following lemma whose proof follows along the lines of Theorem 3.3.3 of [15]. Given $0 \leq m \leq k$, let $\Lambda_m$ be the set of all multiindices $\alpha$ satisfying $|\alpha| \le m$. For a multiindex $\gamma$, denote by $m(\gamma) = \sharp\{\gamma_0 : |\gamma_0| \le |\gamma|\}$. Also for a $|\gamma|$-times differentiable function $\Psi : \mathbb{R}^d \to \mathbb{R}$, denote by $\partial^{\le |\gamma|}\Psi(x)$ the $m(\gamma)$-dimensional vector with entries $\partial^{\gamma_0}\Psi(x)$, $|\gamma_0| \le |\gamma|$. If $\Psi = (\Psi_1, \Psi_2, \ldots, \Psi_d) : \mathbb{R}^d \to \mathbb{R}^d$ is such that each $\Psi_i$ is $|\gamma|$-times continuously differentiable then $\partial^{\le |\gamma|}\Psi(x) \doteq (\partial^{\le |\gamma|}\Psi_1(x), \ldots, \partial^{\le |\gamma|}\Psi_d(x))$. We will call a map $P : \mathbb{R}^m \to \mathbb{R}^d$ a polynomial of degree at most $\wp$ if $P(x) = (P_1(x), \ldots, P_d(x))'$ and each $P_i : \mathbb{R}^m \to \mathbb{R}$ is a polynomial of degree at most $\wp$. Also for $u, v \in \mathbb{R}^l$ we define

$$u * v \doteq (u_1 v_1, \ldots, u_l v_l)'.$$

**Lemma 4.8.** *Let $\alpha, \beta, \gamma$ be multiindices such that $|\alpha|, |\beta|, |\gamma| \le k$. Then there exist subsets $\Lambda^1_\alpha, \Lambda^2_\alpha$ of $\Lambda_{|\alpha|}$ and $\Lambda_{|\alpha|-1}$, respectively, a subset $\Gamma^\alpha_{\beta,\gamma}$ of $\Lambda_{|\gamma|}$ and polynomials $P^\alpha_{\beta,\gamma} : \mathbb{R}^{m(\gamma)} \to \mathbb{R}^d$ of degree at most $|\alpha|$, such that $\partial^\alpha \phi^n$ satisfies:*

$$\begin{aligned}
\partial^\alpha \phi^n_t(x) &= \partial^\alpha x + \int_0^t G^n(\partial^\alpha \phi^n_r(x), \phi^n_r(x), \mathrm{d}r) \\
&\quad + \sum_{(\beta,\gamma) \in \Lambda^1_\alpha \times \Lambda^2_\alpha} \int_0^t G^{\alpha,n}_{\beta,\gamma}(\partial^{\le |\gamma|}\phi^n_r(x), \phi^n_r(x), \mathrm{d}r),
\end{aligned} \qquad (4.11)$$

*where for $x, y \in \mathbb{R}^d, G^n(x, y, r) = \nabla_y \hat{F}^n(y, r) \cdot x$ and for $(x, y) \in \mathbb{R}^{m(\gamma)} \times \mathbb{R}^d, G^{\alpha,n}_{\beta,\gamma}(x, y, r) = P^\alpha_{\beta,\gamma}(x) * \partial^\beta_y \hat{F}^n(y, r)$.*

Note in particular that in the third term on the right-hand side of (4.11), one finds partial derivatives of $\phi^n_r(x)$ of order strictly less than $|\alpha|$.

**Lemma 4.9.** *For each $p \ge 1, L \in (0, \infty)$, there is a constant $k_1 \equiv k_1(k, p, L) \in (0, \infty)$ such that for every multiindex $\alpha, |\alpha| \le k$,*

$$\sup_n \sup_{|x| \le L} \mathbb{E} \sup_{0 \le t \le T} |\partial^\alpha \phi^n_t(x)|^p \le k_1 \qquad (4.12)$$

$$\sup_n \sup_{|x| \le L} \mathbb{E}|\partial^\alpha \phi^n_t(x) - \partial^\alpha \phi^n_{t'}(x)|^p \le k_1 |t - t'|^{p/2}. \qquad (4.13)$$

**Proof.** Fix $L > 0$ and consider $x \in \mathbb{R}^d$ such that $|x| \le L$. We will first show inequality (4.12). It suffices to prove (4.12) for $\alpha = 0$ and establish that if, for some $m < k$, it holds for $\partial^\alpha \phi^n_t$ with $|\alpha| \le m$ and all $p \ge 1$ then it also holds for $\partial_i \partial^\alpha \phi^n_t$ with all $p \ge 1$ (with a possibly larger constant $k_1$) and $i = 1, \ldots, d$. The desired result then follows by induction.

Consider first $\alpha = 0$. For this case the bound in (4.12) follows immediately on using (3.4) and applying the Bürkholder–Davis–Gundy inequality to the square-integrable martingale $N_t = \int_0^t M(\phi^n_r(x), \mathrm{d}r)$ [note that $\langle\!\langle N \rangle\!\rangle_t = \int_0^t a(\phi^n_r(x), \phi^n_r(x), r) \mathrm{d}r$ and $a \in \widetilde{\mathbb{C}}^{k,\delta}_T(\mathbb{R}^{d \times d})$].

Now, suppose that (4.12) holds for all multiindices $\alpha$ with $|\alpha| \leq m$, for some $m < k$. Fix $\alpha$ with $|\alpha| \leq m$, an $i \in \{1, 2, \ldots, d\}$, and consider the multiindex $\tilde{\alpha} = \alpha + 1_i$, where $1_i$ is a $d$-dimensional vector with 1 in the $i$th entry and 0 elsewhere. From Lemma 4.8, one finds that $\partial^{\tilde{\alpha}} \phi^n_t$ solves (4.11) for $\alpha = \tilde{\alpha}$. Note that for $\beta \in \Lambda^1_{\tilde{\alpha}}$,

$$\partial^\beta_y \hat{F}^n(y, t) = \int_0^t \partial^\beta_y b_u(y, s) \, \mathrm{d}s + \sqrt{\varepsilon_n} \partial^\beta_y M(y, t).$$

From (3.4) and recalling that $(b, a) \in \mathbb{C}^{k,\delta}_T(\mathbb{R}^d) \times \widetilde{\mathbb{C}}^{k,\delta}_T(\mathbb{R}^{d \times d})$, we have that for some $c_1, c_2 \in (0, \infty)$,

$$\sup_{0 \leq t \leq T} \sup_{y \in \mathbb{R}^d} \left| \int_0^t \partial^\beta_y b_u(y, s) \, \mathrm{d}s \right| \leq c_1 \quad \text{and} \quad \sup_{0 \leq t \leq T} \sup_{y \in \mathbb{R}^d} |\langle\!\langle \partial^\beta_y M(y, t) \rangle\!\rangle_t| \leq c_2.$$

This along with the assumption

$$\sup_n \sup_{|x| \leq L} \mathbb{E} \sup_{0 \leq t \leq T} |\partial^\nu \phi^n_t(x)|^p \leq k_1 \quad \text{for } \nu, |\nu| \leq |\alpha|,$$

shows that for some $c_3 \in (0, \infty)$, for all $(\beta, \gamma) \in \Lambda^1_{\tilde{\alpha}} \times \Lambda^2_{\tilde{\alpha}}$,

$$\sup_n \sup_{|x| \leq L} \mathbb{E} \sup_{0 \leq t \leq T} \left| \int_0^t G^{\tilde{\alpha}, n}_{\beta, \gamma}(\partial^{\leq |\gamma|} \phi^n_r(x), \phi^n_r(x), \mathrm{d}r) \right|^p \leq c_3.$$

In a similar manner one has for some $c_4 \in (0, \infty)$,

$$\mathbb{E} \sup_{0 \leq t \leq T} \left| \int_0^s G^n(\partial^{\tilde{\alpha}} \phi^n_r(x), \phi^n_r(x), \mathrm{d}r) \right|^p \leq c_4 \int_0^t \mathbb{E}\left( \sup_{0 \leq r \leq s} |\partial^{\tilde{\alpha}} \phi^n_r(x)|^p \right) \mathrm{d}s.$$

Combining the above inequalities we obtain

$$\sup_n \sup_{|x| \leq L} \mathbb{E} \sup_{0 \leq s \leq t} |\partial^{\tilde{\alpha}} \phi^n_s(x)|^p \leq c_3 + c_4 \sup_n \sup_{|x| \leq L} \int_0^t \mathbb{E}\left( \sup_{0 \leq r \leq s} |\partial^{\tilde{\alpha}} \phi^n_r(x)|^p \right) \mathrm{d}s.$$

Now an application of Gronwall's lemma shows that for some $c_5 \in (0, \infty)$,

$$\sup_n \sup_{|x| \leq L} \mathbb{E} \sup_{0 \leq t \leq T} |\partial^{\tilde{\alpha}} \phi^n_t(x)|^p \leq c_5.$$

This establishes (4.12) for all $\tilde{\alpha}$ with $|\tilde{\alpha}| \leq |\alpha| + 1$. Finally consider (4.13). For $t, t' \in [0, T]$, $t' \leq t$, we have from (4.11) that

$$\begin{aligned}
\partial^\alpha \phi^n_t(x) - \partial^\alpha \phi^n_{t'}(x) &= \int_{t'}^t G^n(\partial^\alpha \phi^n_r(x), \phi^n_r(x), \mathrm{d}r) \\
&\quad + \sum_{(\beta, \gamma) \in \Lambda^1_\alpha \times \Lambda^2_\alpha} \int_{t'}^t G^{\alpha, n}_{\beta, \gamma}(\partial^{\leq |\gamma|} \phi^n_r(x), \phi^n_r(x), \mathrm{d}r).
\end{aligned} \tag{4.14}$$

Using (4.12) on the right-hand side of (4.14) we now have (4.13) via an application of Hölder's and Bürkholder–Davis–Gundy's inequalities. □

The proof of Theorem 3.5 proceeds along the lines of Section 5.4 of [15]. We begin by introducing certain Sobolev spaces. Let $j$ be a non-negative integer and let $1 < p < \infty$. Let $B_N \equiv B(0, N)$ be the $\mathbb{R}^d$-ball with center the origin and radius $N$. Let $h \colon \mathbb{R}^d \to \mathbb{R}^d$ be a function such that the distributional derivative $\partial^\alpha h \in L^p(B_N)$ for all $\alpha$ such that $|\alpha| \le j$. Define

$$\|h\|_{j,p:N} = \left(\sum_{|\alpha| \le j} \int_{B_N} |\partial^\alpha h(x)|^p \, dx\right)^{1/p}.$$

The space $H_{j,p}^{\mathrm{loc}} = \{h \colon \mathbb{R}^d \to \mathbb{R}^d, \|h\|_{j,p:N} < \infty \text{ for all } N\}$ together with the seminorms defined above is a real separable semireflexive Fréchet space. By Sobolev's imbedding theorem, we have $H_{j+1,p}^{\mathrm{loc}} \subset \mathbb{C}^j(\mathbb{R}^d) \subset H_{j,p}^{\mathrm{loc}}$ if $p > d$. Furthermore the imbedding $i \colon H_{j+1,p}^{\mathrm{loc}} \to \mathbb{C}^j(\mathbb{R}^d)$ is a compact operator by the Rellich–Kondrachov theorem (see [1]).

**Proposition 4.10.** *The sequence $\{(\phi^n, \hat{F}^n)\}_{n \ge 1}$ is tight in $W_{k-1} \times W_{k-1}$*

**Proof.** It suffices to show that both $\{\phi^n\}_{n \ge 1}$ and $\{\hat{F}^n\}_{n \ge 1}$ are tight in $W_{k-1}$. We will use Kolmogorov's tightness criterion [15], Theorem 1.4.7, page 38. From Lemmas 4.7 and 4.9, we have that for each $p > 1, N > 1$, there exist $c_1, c_2 \in (0, \infty)$ such that for all $t, t' \in [0, T]$

$$\sup_n \mathbb{E}\|\phi_t^n - \phi_{t'}^n\|_{k,p:N}^p \le c_1 |t - t'|^{p/2},$$

$$\sup_n \mathbb{E}\|\hat{F}^n(\cdot, t) - \hat{F}^n(\cdot, t')\|_{k,p:N}^p \le c_2 |t - t'|^{p/2}.$$

Furthermore, since $\hat{F}^n(\cdot, 0) = 0$ and $\phi_0^n(x) = x$, we get that for each $p > 1, N > 1$ there exist $c_3, c_4 \in (0, \infty)$ such that

$$\sup_n \mathbb{E}\|\phi_t^n\|_{k,p:N}^p \le c_3 \quad \text{and} \quad \sup_n \mathbb{E}\|\hat{F}^n(t)\|_{k,p:N}^p \le c_4.$$

Applying Theorem 1.4.7 of [15] with $p > 2$ now gives a tightness of $\{\phi^n\}_{n \ge 1}$ and $\{\hat{F}^n\}_{n \ge 1}$ in the semi-weak topology on $C([0, T] : H_{k,p}^{\mathrm{loc}})$ (cf. [15]). Since the imbedding map $i \colon H_{k,p}^{\mathrm{loc}} \to \mathbb{C}^{k-1}$ is compact, tightness in $W_{k-1} \times W_{k-1}$ with the topology introduced in Section 2 follows (see [15], pages 246–247). □

Recall the definitions (3.2) and (3.3). For the proof of the following lemma we refer the reader to Section 2.1 of [4].

**Lemma 4.11.** *Let $f_n, f \in \hat{W}_{k-1}$ be such that $\sup_{0 \le t \le T} \lambda_{k-1}(f_n(t), f(t)) \to 0$, as $n \to \infty$. Then $\sup_{0 \le t \le T} d_{k-1}(f_n(t), f(t)) \to 0$.*

**Proof of Theorem 3.5.** Convergence as $\mathbb{C}^{k-1}$-flows is immediate from Theorem 4.2, Propositions 4.6 and 4.10. Using Skorokhod's representation theorem, one can find a sequence of pairs $\{(\tilde{\phi}^n, \tilde{F}^n)\}_{n\geq 1}$ that has the same distribution as $\{(\phi^n, \hat{F}^n)\}_{n\geq 1}$ and $\{(\tilde{\phi}^0, \tilde{F}^0)\}$ that has the same distribution as $\{(\phi^0, F^0)\}$ and $\sup_{0\leq t\leq T}[\lambda_k(\tilde{\phi}_t^n, \tilde{\phi}_t^0) + \lambda_k(\tilde{F}^n(t), \tilde{F}^0(t))] \to 0$, a.s. Since $\phi^n, \phi^0 \in \hat{W}_k$ a.s., the same holds for $\tilde{\phi}^n, \tilde{\phi}^0$. Thus from Lemma 4.11 $\sup_{0\leq t\leq T} d_{k-1}(\tilde{\phi}_t^n, \tilde{\phi}_t^0) \to 0$ a.s. Hence $(\phi^n, \hat{F}^n) \to (\phi^0, F^0)$ as $G^{k-1}$-flows. $\square$

## 5. Application to image analysis

A common approach to image matching problems (see [12, 14, 17] and references therein) is to consider a $\mathbb{R}^p$-valued, continuous and bounded function $T(\cdot)$, referred to as the "template" function, defined on a bounded open set $\mathcal{O} \subseteq \mathbb{R}^3$, which represents some canonical example of a structure of interest. By considering all possible smooth transformations $h: \mathcal{O} \to \mathcal{O}$ one can generate a rich library of targets (or images) given by the form $T(h(\cdot))$.

In typical situations we are given data generated by an a priori unknown function $h$, and the key question of image matching is that of estimating $h$ from the observed data. A Bayesian approach to this problem requires a prior distribution on the space of transformations and a formulation of a noise/data model. The "maximum" of the posterior distribution on the space of transformations given the data can then be used as an estimate $\hat{h}$ for the underlying unknown transformation $h$. In certain applications (e.g., medical diagnosis), the goal is to obtain numerical approximations for certain key structures present in the image, such as volumes of subregions, curvatures and surface areas. If the prior distribution on the transformations (and in particular the estimated transformation) is on the space of diffeomorphisms, then this information can be recovered from the template. Motivated by such a Bayesian approach, a variational problem on the space of $\mathbb{C}^m$-diffeomorphic flows was formulated and analyzed in [12].

Before going into the description of this variational problem, we note that although the chief motivation for the variation problem studied in [12] came from Bayesian considerations, no rigorous results on relationships between the two formulations (variational and Bayesian) were established. The goal of our study is to develop a rigorous asymptotic theory that connects a Bayesian formulation for such an image matching problem with the variational approach taken in [12]. The precise result that we will establish is Theorem 5.1, given at the end of this section. The result is an application of Theorem 3.2 for local characteristics $(a, b) \equiv (a, 0)$ and $a \in \widetilde{\mathbb{C}}_T^{k,1/2}(\mathbb{R}^{3\times 3})$, with $k = m - 2$.

Let $\mathcal{C}_0^\infty(\mathcal{O})$ be the space of infinitely differentiable, real-valued functions on $\mathcal{O}$ with compact support in $\mathcal{O}$. The starting point of the variational formulation is a differential operator $L$ on $[\mathcal{C}_0^\infty(\mathcal{O})]^3$, the exact form of which is determined from specific features of the problem under study. The formulation, particularly for problems from biology, often uses principles from physics and continuum mechanics as a guide in the selection of $L$. We refer the reader to Christensen *et al.* [9, 10], where natural choices of $L$ in shape models from anatomy are provided.

Define the norm $\|\cdot\|_L$ on $[\mathcal{C}_0^\infty(\mathcal{O})]^3$ by

$$\|f\|_L^2 \doteq \sum_{i=1}^3 \int_\mathcal{O} |(Lf)_i(u)|^2 \, du,$$

where we write a function $g \in [\mathcal{C}_0^\infty(\mathcal{O})]^3$ as $(g_1, g_2, g_3)'$. It is assumed that $\|\cdot\|_L$ generates an inner product on $[\mathcal{C}_0^\infty(\mathcal{O})]^3$ and that the Hilbert space $H$ defined as the closure of $[\mathcal{C}_0^\infty(\mathcal{O})]^3$ with this inner product is separable. We will need the functions in $H$ to have sufficient regularity and thus assume that the norm $\|\cdot\|_L$ dominates an appropriate Sobolev norm. More precisely, let $W_0^{m+2,2}(\mathcal{O})$ be the closure of $\mathcal{C}_0^\infty(\mathcal{O})$ with respect to the norm

$$\|g\|_{W_0^{m+2,2}(\mathcal{O})} \doteq \left( \int_\mathcal{O} \sum_{|\alpha| \leq m+2} |\partial^\alpha g(u)|^2 \, du \right)^{1/2}, \qquad g \in \mathcal{C}_0^\infty(\mathcal{O}), \tag{5.1}$$

where $\alpha$ denotes a multiindex and $m \geq 3$. Define $\mathcal{V}_m \doteq [W_0^{m+2,2}(\mathcal{O})]^{\otimes 3}$, where $\otimes$ is used to denote the usual tensor product of Hilbert spaces. We denote by $\|\cdot\|_{\mathcal{V}_m}$ the norm on $\mathcal{V}_m$. The main regularity condition on $L$ is the following domination requirement on the $\|\cdot\|_L$ norm. There exists a constant $c \in (0, \infty)$ such that

$$\|f\|_L \geq c \|f\|_{\mathcal{V}_m} \qquad \text{for all } f \in [\mathcal{C}_0^\infty(\mathcal{O})]^3.$$

This condition ensures that $H \subseteq \mathbb{C}^{m,1/2}(\overline{\mathcal{O}})$ (see [1], Theorem 4.12, parts II and III, page 85). We denote by $\mathcal{H}$ the Hilbert space $L^2([0,1]:H)$. For a fixed $\vartheta \in \mathcal{H}$ let $\{\eta_{s,t}(x)\}_{s \leq t \leq 1}$ be the unique solution of the ordinary differential equation

$$\frac{\partial \eta_{s,t}(x)}{\partial t} \doteq \vartheta(\eta_{s,t}(x), t), \qquad \eta_{s,s}(x) = x, \qquad 0 \leq s \leq t \leq 1. \tag{5.2}$$

Then it follows that $\{\eta_{s,t}, 0 \leq s \leq t \leq T\}$ is a forward flow of $\mathbb{C}^m$-diffeomorphisms on $\mathcal{O}$ (see [15], Theorem 4.6.5, page 173). Since $\vartheta(\cdot, t)$ has a compact support in $\overline{\mathcal{O}}$, one can extend $\eta_{s,t}$ to all of $\mathbb{R}^3$ by setting $\eta_{s,t}(x) \equiv x$, if $x \in \mathcal{O}^c$. Extended in this way $\eta_{s,t}$ can be considered as an element of $G^m$, as defined in Section 3. Denoting $\eta_{0,1}$ by $h_\vartheta$, we can now generate a family of smooth transformations (diffeomorphisms) on $\mathcal{O}$ by varying $\vartheta \in \mathcal{H}$. Specifically, the library of transformations that is used in the variational formulation of the image matching problem is $\{h_\vartheta \mid \vartheta \in \mathcal{H}\}$.

We now describe the data that is used in selecting the transformation $h_{\vartheta^*}$ for which the image $T(h_{\vartheta^*}(\cdot))$ best matches the data. Let $\mathcal{L}$ be a finite index set and $\{\mathbb{X}_i\}_{i \in \mathcal{L}}$ be a collection of disjoint subsets of $\mathcal{O}$ such that $\bigcup_{i \in \mathcal{L}} \mathbb{X}_i = \mathcal{O}$. Collected data $\{d_i\}_{i \in \mathcal{L}}$ represents integrated responses over each of the subsets $\mathbb{X}_i, i \in \mathcal{L}$. More precisely, if $T(h(\cdot))$ was the true underlying image and the data were completely error free and noiseless, then $d_i = \int_{\mathbb{X}_i} T(h(\sigma)) \, d\sigma / \mathrm{vol}(\mathbb{X}_i)$, $i \in \mathcal{L}$, where vol denotes volume. Let $d = (d_1, d_2, \ldots, d_n)'$, where $n = |\mathcal{L}|$. Defining $Y_d(x) = d_i$,

$x \in \mathbb{X}_i$, $i \in \mathcal{L}$, the expression

$$\frac{1}{2} \int_{\mathcal{O}} |T(h_\vartheta(x)) - Y_d(x)|^2 \, dx$$

is a measure of discrepancy between a candidate target image $T(h_\vartheta(\cdot))$ and the observations. This suggests a natural variational criterion for selecting the "best" transformation matching the data. The objective function that is minimized in the variational formulation of the image matching problem is a sum of two terms, the first reflecting the "likelihood" of the transformation or change-of-variable $h_\vartheta$ and the second measuring the conformity of the transformed template with the observed data. More precisely, define for $\vartheta \in \mathcal{H}$,

$$J_d(\vartheta) \doteq \frac{1}{2}\left( \|\vartheta\|_\mathcal{H}^2 + \int_{\mathcal{O}} |T(h_\vartheta(x)) - Y_d(x)|^2 \, dx \right). \tag{5.3}$$

Then $\vartheta^* \in \operatorname{argmin}_{\vartheta \in \mathcal{H}} J_d(\vartheta)$ represents the "optimal" velocity field that matches the data $d$ and for which the $h_{\vartheta^*}$, obtained by solving (5.2), gives the "optimal" transformation. This transformation then yields an estimate of the target image as $T(h_{\vartheta^*}(\cdot))$. Equivalently, defining for each $h \in G^0$,

$$\hat{J}_d(h) \doteq \inf_{\vartheta \in \Psi_h} J_d(\vartheta) \qquad (\text{where } \Psi_h = \{\vartheta \in \mathcal{H} : h = h_\vartheta\}), \tag{5.4}$$

we see that an optimal transformation is $h^* = h_{\vartheta^*} \in \operatorname{argmin}_h \hat{J}_d(h)$.

Up to a relabeling of the time variable, the above variational formulation (in particular the cost function in (5.3)) was motivated in [12] through Bayesian considerations, but no rigorous justification was provided. In [12] the orientation of time is consistent with the change-of-variable evolving toward the identity mapping at the terminal time. To relate the variational problem to stochastic flows it is more convenient to have the identity mapping at time zero. We next introduce a stochastic Bayesian formulation of the image matching problem and describe the precise asymptotic result that we will establish.

Let $\{\phi_i\}$ be a complete orthonormal system in $H$ and $\beta \equiv (\beta_i)_{i=1}^\infty$ be as in Section 2, a sequence of independent, standard, real-valued Brownian motions on some filtered probability space $(\Omega, \mathcal{F}, \mathbb{P}, \{\mathcal{F}_t\})$. Write $(C([0,T] : \mathbb{R}^\infty), \mathcal{B}(C([0,T] : \mathbb{R}^\infty))) \equiv (S, \mathcal{S})$, and note that $\beta$ is a random variable with values in $S$. Consider the stochastic flow

$$d\psi_{s,t}(x) = \sqrt{\varepsilon} \sum_{i=1}^\infty \phi_i(\psi_{s,t}(x)) \, d\beta_i(t), \qquad \psi_{s,s}(x) = x, \qquad x \in \mathcal{O}, 0 \le s \le t \le 1, \tag{5.5}$$

where $\varepsilon \in (0, \infty)$ is fixed. From Maurin's theorem (see [1], Theorem 6.61, page 202) it follows that the imbedding map $H \to \mathcal{V}_{m-2}$ is Hilbert–Schmidt. Also, $\mathcal{V}_{m-2}$ is continuously embedded in $\mathbb{C}^{m-2,1/2}(\overline{\mathcal{O}})$. Thus for some $k_1, k_2 \in (0, \infty)$ and all $u, x, y \in \mathcal{O}$,

$$\sum_{i=1}^\infty |\phi_i(u)|^2 \le k_1 \sum_{i=1}^\infty \|\phi_i\|_{\mathcal{V}_{m-2}}^2 < \infty,$$

$$\sum_{i=1}^{\infty} |\phi_i(x) - \phi_i(y)|^2 \leq k_1 |x-y|^2 \sum_{i=1}^{\infty} \|\phi_i\|_{\mathcal{V}_{m-2}}^2 = k_2 |x-y|^2.$$

One also has that if $\phi_l$ is extended to all of $\mathbb{R}^3$ by setting $\phi_l(u) = 0$, for all $x \in \mathcal{O}^c$, then $a(x,y) = \sum_{l=1}^{\infty} \phi_l(x) \phi_l'(y)$ is in $\tilde{\mathbb{C}}_T^{m-2,1/2}(\mathbb{R}^{3\times 3})$. Thus it follows (cf. [15], pages 80 and 106) that

$$F(x,t) = \sum_{l=1}^{\infty} \int_0^t \phi_i(x) \, \mathrm{d}\beta_i(r)$$

is a $\mathbb{C}^{m-2,\nu}$-Brownian motion, $0 < \nu < 1/2$, with local characteristics $(a,0)$. Also (5.5) admits a unique solution $\{\psi_{s,t}^\varepsilon(x), 0 \leq s \leq t \leq 1\}$ for each $x \in \mathcal{O}$ and $\{\psi_{s,t}^\varepsilon\}_{0 \leq s \leq t \leq 1}$ is a forward flow of $\mathbb{C}^k$-diffeomorphisms, with $k = m - 2$ (see [15], Theorem 4.6.5). In particular, $X^\varepsilon \doteq \psi_{0,1}^\varepsilon$ is a random variable in the space of $\mathbb{C}^k$-diffeomorpshisms on $\mathcal{O}$. The law of $X^\varepsilon$ (for a fixed $\varepsilon > 0$) on $G^k$ will be used as the prior distribution on the transformation space $G^k$. Note that $T(X^\varepsilon(\cdot))$ induces a measure on the space of target images.

We next consider the data model. Let $\mathcal{L}$ and $n$ be as introduced below (5.2). We suppose that the data is given through an additive Gaussian noise model:

$$D_i = \int_{\mathbb{X}_i} T(X^\varepsilon(x)) \, \mathrm{d}x + \sqrt{\varepsilon} \xi_i,$$

where $\{\xi_i, i \in \mathcal{L}\}$ is a family of independent, $p$-dimensional standard normal random variables.

In the Bayesian approach to the image matching problem one considers the posterior distribution of $X^\varepsilon$ given the data $D$ and uses the "mode" of this distribution as an estimate for the underlying true transformation. More precisely, let $\{\Gamma^\varepsilon\}_{\varepsilon > 0}$ be a family of measurable maps from $\mathbb{R}^{np}$ to $\mathcal{P}(G^k)$ (the space of probability measures on $G^k$), such that

$$\Gamma^\varepsilon(A|D) = \mathbb{P}[X^\varepsilon \in A|D] \qquad \text{a.s. for all } A \in \mathcal{B}(G^k).$$

We refer to $\Gamma^\varepsilon(\cdot|d)$ as a regular conditional probability distribution (r.c.p.d.) of $X^\varepsilon$ given $D = d$. In Theorem 5.1 below, we will show that there is a r.c.p.d. $\{\Gamma^\varepsilon(\cdot|d), d \in \mathbb{R}^{np}\}_{\varepsilon > 0}$ such that for each $d \in \mathbb{R}^{np}$, the family $\{\Gamma^\varepsilon(\cdot|d)\}_{\varepsilon > 0}$, regarded as elements of $\mathcal{P}(G^{k-1}) \supseteq \mathcal{P}(G^k)$, satisfies an LDP with rate function

$$I_d(h) = \hat{J}_d(h) - \lambda_d, \qquad \text{where } \lambda_d = \inf_{h \in G^{k-1}} \hat{J}_d(h) = \inf_{\vartheta \in \mathcal{H}} J_d(\vartheta).$$

Formally writing $\Gamma^\varepsilon(A|d) \approx \int_A \mathrm{e}^{-I_d(h)/\varepsilon} \, \mathrm{d}h$, one sees that for small $\varepsilon$, the "mode" of the posterior distribution given $D = d$, which represents the "optimal transformation" in the Bayesian formulation, can be formally interpreted as $\operatorname{argmin}_h I_d(h)$. Note that $\hat{J}_d(h) = \infty$ if $h \notin G^m$ (recall $m = k + 2$). Theorem 5.1 in particular says that $h \in G^m$ is a $\delta$-minimizer for $I_d(h)$ if and only if it is also a $\delta$-minimizer for $\hat{J}_d(h)$. Thus Theorem 5.1 makes precise the asymptotic relationship between the variational and the Bayesian formulation of the above image matching problem.

**Theorem 5.1.** *There exists an r.c.p.d. $\Gamma^\varepsilon$ such that for each $d \in \mathbb{R}^n$, the family of probability measures $\{\Gamma^\varepsilon(d)\}_{\varepsilon>0}$ on $G^{k-1}$ satisfies a large deviation principle (as $\varepsilon \to 0$) with rate function*

$$I_d(h) \doteq \hat{J}_d(h) - \lambda_d. \tag{5.6}$$

We begin with the following proposition. Let $\tilde{I}: G^{k-1} \to [0, \infty]$ be defined as

$$\tilde{I}(h) \doteq \inf_{\vartheta \in \Psi_h} \tfrac{1}{2}\|\vartheta\|_{\mathcal{H}}^2.$$

**Proposition 5.2.** *The family $\{X^\varepsilon\}_{\varepsilon>0}$ satisfies a LDP in $G^{k-1}$ with rate function $\tilde{I}$.*

Proposition 5.2 is consistent with results in Section 3 in that although the local characteristics are in $\mathbb{C}^k$ and $X^\varepsilon \in G^k$, the LDP is established in the larger space $G^{k-1}$. This is due to the tightness issues described in the Introduction. Furthermore, as noted below (5.2), if $\|\vartheta\|_{\mathcal{H}} < \infty$ then $\vartheta$ induces a flow of $\mathbb{C}^m$-diffeomorphisms on $\mathcal{O}$. Thus if $h \in G^{k-1} \setminus G^m$ then $\Psi_h$ is empty, and consequently $\tilde{I}(h) = \infty$. Hence there is a further widening of the "gap" between the regularity needed for the rate function to be finite and the regularity associated with the space in which the LDP is set. This is due to the fact that the variational problem is formulated essentially in terms of $L^2$ norms of derivatives, while in the theory of stochastic flows as developed in [15] assumptions are phrased in terms of $L^\infty$ norms.

**Proof of Proposition 5.2.** From Theorem 3.2 and an application of the contraction principle we have that $\{X^\varepsilon\}_{\varepsilon>0}$ satisfies LDP in $G^{k-1}$ with rate function

$$I^*(h) \doteq \inf_{u \in \mathcal{L}^*(h)} \frac{1}{2} \int_0^T \|u(s)\|_{l_2}^2 \, \mathrm{d}s,$$

where $\mathcal{L}^*(h) = \{u \in L^2([0,1]:l_2) | h = \phi^{0,u}(1)\}$ and where $\phi^{0,u}$ is defined via (3.6), but with $f_i$ there replaced by $\phi_i$. Note that there is a one-to-one correspondence between $u \in L^2([0,1]:l_2)$ and $\vartheta \in \mathcal{H}$ given as $\vartheta(t,x) = \sum_{l=1}^\infty u_l(t)\phi_l(x)$ and $\int_0^T \|u(s)\|_{l_2}^2 \, \mathrm{d}s = \|\vartheta\|_{\mathcal{H}}^2$. In particular $u \in \mathcal{L}^*(h)$ if and only if $\vartheta \in \Psi_h$. Thus $I^*(h) = \tilde{I}(h)$ and the result follows. □

**Proposition 5.3.** *For each $d \in \mathbb{R}^n$, $I_d$ defined in (5.6) is a rate function on $G^{k-1}$.*

**Proof.** From (5.6) and the definition of $\tilde{I}$ we have for $h \in G^{k-1}$ that

$$I_d(h) = \tilde{I}(h) + \frac{1}{2} \int_{\mathcal{O}} |T(h(x)) - Y_d(x)|^2 \, \mathrm{d}x$$

$$- \inf_{h \in G^{k-1}} \left\{ \tilde{I}(h) + \frac{1}{2} \int_{\mathcal{O}} |T(h(x)) - Y_d(x)|^2 \, \mathrm{d}x \right\}.$$

From Proposition 5.2, $\tilde{I}$ is a rate function and therefore has compact level sets. Additionally $T$ is a continuous and bounded function on $\mathcal{O}$. The result follows. □

**Proof of Theorem 5.1.** We begin by noting that $\Gamma^\varepsilon(\cdot|d)$ defined as

$$\Gamma^\varepsilon(A|d) \doteq \frac{\int_A e^{-1/(2\varepsilon)\sum_{i=1}^n |d_i - \int_{\mathbb{X}_i} T(h(y))\,\mathrm{d}y|^2}\mu^\varepsilon(\mathrm{d}h)}{\int_{G^{k-1}} e^{-1/(2\varepsilon)\sum_{i=1}^n |d_i - \int_{\mathbb{X}_i} T(h(y))\,\mathrm{d}y|^2}\mu^\varepsilon(\mathrm{d}h)},$$

where $\mu^\varepsilon = \mathbb{P} \circ (X^\varepsilon)^{-1} \in \mathcal{P}(G^{k-1})$ is an r.c.p.d. of $X^\varepsilon$ given $D = d$. Using the equivalence between the Laplace principle and large deviations principle (see [11], Section 1.2), it suffices to show that for all continuous and bounded real functions $F$ on $G^{k-1}$,

$$-\varepsilon \log \int_{G^{k-1}} \exp\left[-\frac{1}{\varepsilon}F(v)\right]\Gamma^\varepsilon(\mathrm{d}v|d) \tag{5.7}$$

converges to $\inf_{h \in G^{k-1}}\{F(h) + I_d(h)\}$. Note that (5.7) can be expressed as

$$\begin{aligned}
&-\varepsilon \log \int_{G^{k-1}} e^{-1/\varepsilon[F(h) + 1/2\sum_{i=1}^n |d_i - \int_{\mathbb{X}_i} T(h(y))\,\mathrm{d}y|^2]}\mu^\varepsilon(\mathrm{d}h) \\
&+ \varepsilon \log \int_{G^{k-1}} e^{-1/\varepsilon[1/2\sum_{i=1}^n |d_i - \int_{\mathbb{X}_i} T(h(y))\,\mathrm{d}y|^2]}\mu^\varepsilon(\mathrm{d}h).
\end{aligned} \tag{5.8}$$

From Proposition 5.2 we see that the first term converges to

$$\begin{aligned}
&\inf_{h \in G^{k-1}}\left\{\tilde{I}(h) + F(h) + \frac{1}{2}\sum_{i=1}^n \left|d_i - \int_{\mathbb{X}_i} T(h(y))\,\mathrm{d}y\right|^2\right\} \\
&= \inf_{h \in G^{k-1}} \inf_{\vartheta \in \Psi_h}\left\{F(h) + \frac{1}{2}\|\vartheta\|_{\mathcal{H}}^2 + \frac{1}{2}\int_{\mathcal{O}} |T(h(y)) - Y_d(y)|^2\,\mathrm{d}y\right\} \\
&= \inf_{h \in G^{k-1}}\{F(h) + \hat{J}_d(h)\},
\end{aligned}$$

where the last equality is a consequence of (5.3) and (5.4). Taking $F = 0$ in the above display, we see that the second term in (5.8) converges to $-\lambda_d$. This proves the result. □

## Notational conventions

- Transpose of a $d$-dimensional vector $v$ will be denoted by $v'$.
- $K \subset\subset \mathbb{R}^d$ means that $K$ is a compact subset of $\mathbb{R}^d$.
- Given $k$- and $m$-dimensional continuous local martingales $M, N$ on some filtered probability space $(\Omega, \mathcal{F}, \mathbb{P}, \{\mathcal{F}_t\})$, we write the cross-quadratic variation of $M$ and $N$ as $\langle\!\langle M, N \rangle\!\rangle_t$, and write $\langle\!\langle N \rangle\!\rangle_t$ when $M = N$. This is a continuous $\mathbb{R}^{m \times k}$-valued $\{\mathcal{F}_t\}$-adapted process.

- Borel $\sigma$-fields on a Polish space $\mathcal{E}$ will be written as $\mathcal{B}(\mathcal{E})$.
- Generic constants will be denoted as $c_1, c_2, \ldots$. Their values may change from one proof to the next.
- The infimum over an empty set is taken to be $\infty$.

## Acknowledgements

Amarjit Budhiraja and Vasileios Maroulas' research was supported in part by the Army Research Office (Grants W911NF-04-1-0230, W911NF-0-1-0080). Paul Dupuis' research was supported in part by the National Science Foundation (NSF-DMS-0404806 and NSF-DMS-0706003) and the Army Research Office (W911NF-05-1-0289).